\documentclass[12pt,a4paper]{article}

\usepackage{amsmath}
\usepackage{amssymb}
\usepackage{bbold}

\allowdisplaybreaks

\renewcommand{\leq}{\leqslant}

\renewcommand{\geq}{\geqslant}

\makeatletter
\newcommand*{\bigcorr@macro}[2]{\sbox{0}{\mbox{$#1($}}\dimen0=\ht0
                \advance\dimen0 by \dp0
                \multiply\dimen0 by #2 \divide\dimen0 by 100}
\newcommand*{\bigcorr@big}[2]{\mbox{$#1\left#2\bigcorr@macro{#1}{85}\vrule
                   height \dimen0 depth 0pt width 0pt\right.\n@space$}}
\newcommand*{\bigcorr@Big}[2]{\mbox{$#1\left#2\bigcorr@macro{#1}{115}\vrule
                   height \dimen0 depth 0pt width 0pt\right.\n@space$}}
\newcommand*{\bigcorr@bigg}[2]{\mbox{$#1\left#2\bigcorr@macro{#1}{145}\vrule
                   height \dimen0 depth 0pt width 0pt\right.\n@space$}}
\newcommand*{\bigcorr@Bigg}[2]{\mbox{$#1\left#2\bigcorr@macro{#1}{175}\vrule
                   height \dimen0 depth 0pt width 0pt\right.\n@space$}}
\DeclareRobustCommand*{\big}[1]{{\mathpalette\bigcorr@big{#1}}}
\DeclareRobustCommand*{\Big}[1]{{\mathpalette\bigcorr@Big{#1}}}
\DeclareRobustCommand*{\bigg}[1]{{\mathpalette\bigcorr@bigg{#1}}}
\DeclareRobustCommand*{\Bigg}[1]{{\mathpalette\bigcorr@Bigg{#1}}}

\makeatother

\newtheorem{theorem}{Theorem}
\newtheorem{definition}[theorem]{Definition}
\newtheorem{lemma}[theorem]{Lemma}

\newcommand{\pb}{\mathbf{P}}
\newcommand{\eb}{\mathbf{E}}
\newcommand{\RR}{\mathbb{R}}
\newcommand{\NN}{\mathbb{N}}
\newcommand{\1}{\mathbb{1}}
\newcommand{\limINF}{\mathop{\underline{\lim}}\limits}
\newcommand{\sign}{\mathop{\rm sign}\nolimits}
\newcommand*{\abs}[1]{\left|#1\right|}
\newcommand{\wh}{\widehat}
\newcommand{\wt}{\widetilde}
\newcommand{\adh}{\overline}
\newcommand{\myint}{\int\limits}
\newcommand{\limn}{\lim_{n\rightarrow\infty}}
\newcommand{\phanteq}{\mathrel{\phantom{=}}}

\begin{document}

\title{Estimation of the Location of a\\
       $0$-type or $\infty$-type Singularity\\
       by Poisson Observations}
\author{Sergue\"{\i} \textsc{Dachian}\\
Laboratoire de Math\'ematiques\\
Universit\'e Blaise Pascal\\
63177 Aubi\`ere CEDEX, France\\
Serguei.Dachian@math.univ-bpclermont.fr}

\date{}
\maketitle

\begin{abstract}
We consider an inhomogeneous Poisson process $X$ on $[0{,}T]$.  The intensity
function of $X$ is supposed to be strictly positive and smooth on $[0{,}T]$
except at the point $\theta$, in which it has either a $0$-type singularity
$\bigl($tends to $0$ like $\abs{x}^p$, $p\in(0{,}1)\bigr)$, or an
$\infty$-type singularity $\bigl($tends to~$\infty$ like $\abs{x}^p$,
$p\in(-1{,}0)\bigr)$.  We suppose that we know the shape of the intensity
function, but not the location of the singularity.  We consider the problem of
estimation of this location (shift) parameter~$\theta$ based on $n$
observations of the process $X$.  We study the Bayesian estimators and, in the
case $p>0$, the maximum likelihood estimator.  We show that these estimators
are consistent, their rate of convergence is $n^{1/(p+1)}$, they have
different limit distributions, and the Bayesian estimators are asymptotically
efficient.
\end{abstract}



\section{Introduction}

Inhomogeneous Poisson process is one of the simplest point processes
$\bigl($see, for example, Daley and Vere-Jones~\cite{DV}$\bigr)$. However, due
to the large choice of intensity functions, this model is reach enough and is
widely used in many applied statistical problems, such as optical
communications, reliability, biology, medicine, image treatment, and so on
$\bigl($see, for example, Karr~\cite{Kar1} and~\cite{Kar2}, Snyder and
Miller~\cite{SM} and Thompson~\cite{Th}$\bigr)$.

The diversity of applications is also due to the possibility of using the
likelihood ratio analysis.  In parameter estimation problems the large samples
theory is quite close to the one of the classical (i.i.d.)  statistics. In
particular, let us consider the problem of estimation of the
parameter~$\theta$ by $n$ independent observations on some fixed interval
$[0{,}T]$ of an inhomogeneous Poisson process $X=\{X(t),\ 0\leq t\leq T\}$ of
intensity function $S_\theta(t)$.  Let us mention that this problem is
equivalent to the one of estimation of the parameter by one observation on a
growing interval of a periodic inhomogeneous Poisson process.  If the problem
is regular (the model is locally asymptotically normal), then both the maximum
likelihood estimator (MLE) $\wh\theta_n $ and the Bayesian estimators (BE)
$\wt\theta_n$ are consistent, asymptotically normal:
$$
\sqrt{n}\bigl(\wh\theta_n-\theta\bigr)\Longrightarrow\mathcal{N}
\bigl(0,I(\theta)^{-1}\bigr),\qquad
\sqrt{n}\bigl(\wt\theta_n-\theta\bigr)\Longrightarrow\mathcal{N}
\bigl(0,I(\theta)^{-1}\bigr),
$$
and asymptotically efficient $\bigl($see, for example, Kutoyants~\cite{Kut2}
and~\cite{Kut}$\bigr)$. Here $I(\theta)$ is the Fisher information given by
$$
I(\theta)=\int_{0}^{T}\frac{\dot S^2_\theta(t)}{S_\theta(t)}\;dt,
$$
where $S_\theta(t)$ is the intensity function and $\dot
S_\theta(t)=\frac{\partial}{\partial\theta}S_\theta(t)$.

If the problem is not regular, then the properties of estimators essentially
change. For example, if $S_\theta(\cdot)$ is smooth everywhere on $[0{,}T]$
except at the point $\theta$, in which it has a jump (consider for instance
$S_\theta(t)=s(t-\theta)$ where $s(\cdot)$ is discontinuous in $0$), then the
MLE and BE are still consistent, but converge at a faster rate:
$$
n\bigl(\wh\theta_n-\theta\bigr)\Longrightarrow\xi_1,\qquad
n\bigl(\wt\theta_n-\theta\bigr)\Longrightarrow\xi_2,
$$
have different limit distributions ($\xi_1 $ and $\xi_2 $ are different with
$\eb\xi_1^2 >\eb\xi_2^2$), and the BE are asymptotically efficient
$\bigl($see, for example, Kutoyants~\cite{Kut2} and~\cite{Kut}$\bigr)$.

In this paper we deal with the case where the intensity function
$S_\theta(\cdot)$ is smooth everywhere on $[0{,}T]$ except at the point
$\theta$, in which it has a singularity of order $p$. The cusp type
singularities were already studied in the preceding paper~\cite{Dach}. Here we
consider $0$-type and $\infty$-type singularities. More precisely, we suppose
that $S_\theta(t)=s(t-\theta)$, where $s(\cdot)$ is some known strictly
positive function on $[-T{,}T]\setminus\{0\}$ and $\theta\in(0{,}T)$ is some
unknown parameter, and that we have the following representation
$$
S_\theta(t)=s(t-\theta)=\begin{cases}
\vphantom{|_{\big)}}a\abs{t-\theta}^p+\psi(t-\theta),&\text{if }
t<\theta\\
\vphantom{|^{\big)}}b\abs{t-\theta}^p+\psi(t-\theta),&\text{if }
t>\theta
\end{cases}
\quad,
$$
where $a,b>0$, $p>-1$ (to guarantee the finiteness of intensity measure), and
$\psi(\cdot)$ is smooth.

If $\psi(0)\ne0$ and $p>1/2$ then, in spite of the singularity of the
intensity function in $\theta$, the Fisher information is finite, and so this
case can be treated as the regular one.

If $\psi(0)\ne0$ and $0<p<1/2$ we say that the intensity function has a cusp
at $\theta$. This is the case treated in~\cite{Dach} (where instead of $a,b>0$
it was supposed $a^2+b^2>0$ only). There it was shown that the MLE and the BE
are consistent, converge at the rate $n^{1/(2p+1)}$ (which is faster than in
the regular case but slower than in discontinuous case):
$$
n^{1/(2p+1)}\bigl(\wh\theta_n-\theta\bigr)\Longrightarrow\eta_1,\qquad
n^{1/(2p+1)}\bigl(\wt\theta_n-\theta\bigr)\Longrightarrow\eta_2,
$$
have different limit distributions, and the BE are asymptotically
efficient. The convergence of moments was equally verified.

If $\psi(0)=0$ and $p>1$ then, as above, the Fisher information is finite and
this case can be treated as the regular one.

If $\psi(0)=0$ and $0<p<1$ we say that the intensity function has a $0$-type
singularity at $\theta$. In this case we study the asymptotic behavior of the
MLE and the BE, and we prove that the estimators are consistent, converge at
the rate $n^{1/(p+1)}$ (which is again intermediate between the regular and
discontinuous case rates), have different limit distributions, and the BE are
asymptotically efficient. We verify also the convergence of moments.

If $-1<p<0$ we say that the intensity function has a $\infty$-type singularity
at $\theta$. In this case we study the asymptotic behavior of the BE only (MLE
makes no sense in this case). We prove that the estimators are consistent,
converge at the rate $n^{1/(p+1)}$ (which is even faster than in discontinuous
case), and are asymptotically efficient. We verify as well the convergence of
moments.

Let us note, that the jump can also be considered as a singularity by taking
$p=0$ and $a\ne b$, which explains that the rates are slower for $p>0$ and
faster for $p<0$.

Let us also mention, that our results are similar to those obtained by
Ibragimov and Khasminskii for the problem of estimation of a singularity
location of the density for the i.i.d.\ model of observations. An exhaustive
exposition of the results can be found in Chapter~6 of their book~\cite{IKh},
but one can also refer to their previous works~\cite{IKh2}
and~\cite{IKh3}. The asymptotic behavior of the MLE and of a wide class of BE
obtained for this (i.i.d.)\ model is similar to the one obtained here for the
model of Poisson observations. Particularly, the rate of convergence of the
estimators is $n^{1/(p+1)}$, and the BE are asymptotically efficient.

Finally, let us note that for the study of the asymptotic behavior of the
estimators we use the method of Ibragimov and Khasminskii presented in their
book~\cite{IKh} (see as well Kutoyants~\cite{Kut}, where this method is
applied to inhomogeneous Poisson process).

\section{Main results}

Suppose we observe $n$ realizations $(X_1,\ldots,X_n)=X^n$ of the Poisson
process $X=\{X(t),\ 0\leq t\leq T\}$ of intensity function
$S_\theta(t)=s(t-\theta)$, where $\theta$ is some unknown parameter,
$\theta\in\Theta=(\alpha{,}\beta)\subseteq(0{,}T)$, and $s(\cdot)$ is some
known strictly positive function on $[-T{,}T]\setminus\{0\}$. We suppose that
the function $s(\cdot)$ can be written in the form $s(t)=d(t)\abs
t^p+\psi(t)$, where $p\in(-1{,}0)\cup(0{,}1)$,
$$
d(t)=\begin{cases}
\vphantom{|_{\big)}}a,&\text{if }t<0\\
\vphantom{|^{\big)}}b,&\text{if }t>0
\end{cases}\quad,
$$
$a,b>0$, and the function $\psi(\cdot)$ is H\"older continuous on $[-T{,}T]$ of
order higher than $(p+1)/2\,$, that is $\bigl|\psi(x)-\psi(y)\bigr|\leq
L\,\abs{x-y}^{\varkappa}$ for all $x,y\in[-T,T]$ with some fixed constants
$L>0$ and $\varkappa>(p+1)/2$. In the case $p>0$ we suppose equally that
$\psi(0)=0$. Our aim is to estimate the parameter~$\theta$ and to study the
asymptotic behavior of estimators as $n$ goes to infinity.

The likelihood ratio in our problem can be written $\bigl($see, for
example,~\cite{Kut}$\bigr)$ as
\begin{align*}
L(\theta,\theta_1,X^n)=\exp\Biggl\{&\sum_{i=1}^{n}\myint_{0}^{T}
\ln\frac{S_\theta(t)}{S_{\theta_1}(t)}\;dX_i(t)\\
&-n\myint_{0}^{T}
\left[\frac{S_\theta(t)}{S_{\theta_1}(t)}-1\right]S_{\theta_1}(t)\;dt\Biggr\},
\end{align*}
where $\theta_1$ is some fixed value of $\theta$.

As usually, introduce the MLE $\wh\theta_n$ as one of the solutions of the
equation
$$
L(\wh\theta_n,\theta_1,X^n)=\sup_{\theta\in\Theta}L(\theta,\theta_1,X^n)
$$
and the BE $\wt\theta_n$ for prior density $q$ (supposed to be positive and
continuous) and quadratic loss function as
$$
\wt\theta_n=\myint_{\alpha}^{\beta}\theta\,q\bigl(\theta|X^n\bigr)\;d\theta,
$$
where the posterior density
$$
q\bigl(\theta|X^n\bigr)=L(\theta,\theta_1,X^n)\,q(\theta)\left(
\myint_{\alpha}^{\beta}L(\theta,\theta_1,X^n)\,q(\theta)\;d\theta\right)^{-1}.
$$

Note that the MLE makes no sense in the case $p<0$, since in this case the
likelihood equals infinity in any point $\theta$ which is event of one of the
Poisson processes $X_1,\ldots,X_n$.

To describe the properties of these estimators we need to introduce the
stochastic process
\begin{align*}
Z(u)=\exp\Biggl\{&p\myint_{-\infty}^{+\infty}\ln\abs{1-\frac uz}\;\pi(dz)+\ln
\frac ab \myint_{0}^{u}\;Y(dz)\\
&-\myint_{-\infty}^{+\infty}\biggl[\abs{1-\frac uz}^p-1-p\,\ln\abs{1-\frac uz}
\biggr] d(z)\abs z^p\;dz\\
&-\frac{a-b}{p+1}\,\abs u^{p+1}\,\sign(u)\Biggr\},
\quad u\in\RR.
\end{align*}
Here and in the sequel $Y$ denotes a Poisson process on $\RR$ of intensity
function $S_0(z)=d(z)\abs z^p$, and $\pi$ is its centered version~:
$\pi=Y-\eb Y$.

We introduce also the random variable $\zeta$, and in the case $p>0$ the random
variable $\xi$ by the equations
$$
\zeta=\myint_{-\infty}^{+\infty}u\,Z(u)\;du\;
\left(\:\myint_{-\infty}^{+\infty}Z(u)\;du\right)^{-1}
$$
and
$$
Z(\xi)=\sup_{u\in\RR}Z(u).
$$
Let us note here, that $\xi$ is well defined in the case $p>0$, since in this
case with probability one the process $Z(u)$ attains its maximum in a unique
point $\bigl($see, for example,~\cite{Erm}$\bigr)$.

Now we can finally state the main results of this paper.

\begin{theorem}
\label{T1}
Under the maid assumptions, the following lower bound on the risks of all
estimators holds: for any $\theta_0\in\Theta$ we have
$$
\lim_{\delta\rightarrow0}\ \limINF_{n\rightarrow\infty}\
\inf_{\adh\theta_n}\ \sup_{\left|\theta-\theta_0\right|<\delta}\
\eb_{\theta}\Bigl(n^{1/(p+1)}\bigl(\adh\theta_n-\theta\bigr)\Bigr)^2
\geq\eb\zeta^2,
$$
where $\inf$ is taken over all possible estimators $\adh\theta_n$ of~$\theta$.
\end{theorem}

This theorem leads us to introduce the following
\begin{definition}
We say that the estimator $\adh\theta_n$ is asymptotically efficient if
$$
\lim_{\delta\rightarrow0}\ \limn\
\sup_{\left|\theta-\theta_0\right|<\delta}\
\eb_{\theta}\Bigl(n^{1/(p+1)}\bigl(\adh\theta_n-\theta\bigr)\Bigr)^2
=\eb\zeta^2
$$
for any $\theta_0\in\Theta$.
\end{definition}

For the BE we have the following
\begin{theorem}
\label{T2}
The BE $\wt\theta_n$ have uniformly in $\theta\in\mathbf{K}$ (for any compact
$\mathbf{K}\subset\Theta$) the following properties:
\begin{itemize}
\item $\wt\theta_n$ is consistent, that is
$$
\wt\theta_n\stackrel{\pb_\theta}{\longrightarrow}\theta
\mbox{ (convergence in probability),}
$$
\item the limit distribution of $\wt\theta_n$ is $\zeta$, that is
$$
n^{1/(p+1)}\bigl(\wt\theta_n-\theta\bigr)\Longrightarrow\zeta
\mbox{ (convergence in law),}
$$
\item for any $k>0$ we have
$$
\limn\eb_\theta
\abs{n^{1/(p+1)}\bigl(\wt\theta_n-\theta\bigr)}^k=\eb\abs\zeta^k
$$
and, moreover, $\wt\theta_n$ is asymptotically efficient.
\end{itemize}
\end{theorem}

And for the MLE (in the case $p>0$) we have the following
\begin{theorem}
\label{T3}
Let $p\in(0{,}1)$. The MLE $\wh\theta_n$ has uniformly in
$\theta\in\mathbf{K}$ (for any compact $\mathbf{K}\subset\Theta$) the
following properties:
\begin{itemize}
\item $\wh\theta_n$ is consistent, that is
$$
\wh\theta_n\stackrel{\pb_\theta}{\longrightarrow}\theta
\mbox{ (convergence in probability),}
$$
\item the limit distribution of $\wh\theta_n$ is $\xi$, that is
$$
n^{1/(p+1)}\bigl(\wh\theta_n-\theta\bigr)\Longrightarrow
\xi\mbox{ (convergence in law),}
$$
\item for any $k>0$ we have
$$
\limn\eb_\theta
\abs{n^{1/(p+1)}\bigl(\wh\theta_n-\theta\bigr)}^k=\eb\abs\xi^k.
$$
\end{itemize}
\end{theorem}

To prove the above stated theorems we apply the method of Ibragimov and
Khasminskii $\bigl($see~\cite{IKh}$\bigr)$. For this we denote
$\theta_u=\theta+u\,n^{-1/(p+1)}$ for all $u\in
U_n=\left(n^{1/(p+1)}(\alpha-\theta)\,{,}\,n^{1/(p+1)}(\beta-\theta)\right)$,
we introduce the normalized likelihood ratio process as
$$
Z_n(u)=L\left(\theta_u,\theta,X^n\right),\quad u\in U_n,
$$
and we establish (the proofs are in the next section) the following three
lemmas.

\begin{lemma}
\label{L1}
The finite-dimensional distributions of $Z_n(u)$ converge to those of $Z(u)$
uniformly in $\theta\in\mathbf{K}$ (for any compact
$\mathbf{K}\subset\Theta$).
\end{lemma}

\begin{lemma}
\label{L2}
For any compact $\mathbf{K}\subset\Theta$ there exists some positive
constant~$C$ such that
$$
\eb_\theta
\left|Z_n^{1/2}(u_1)-Z_n^{1/2}(u_2)\right|^2\leq C\,\abs{u_1-u_2}^{p+1}
$$
for all $u_1,u_2\in U_n$, $\theta\in\mathbf{K}$ and $n$ sufficiently large.
\end{lemma}

\begin{lemma}
\label{L3}
For any compact $\mathbf{K}\subset\Theta$ there exists some positive
constant~$c$ such that
$$
\eb_\theta Z_n^{1/2}(u)\leq\exp\bigl\{-c\,\abs u^{p+1}\bigr\}
$$
for all $u\in U_n$, $\theta\in\mathbf{K}$ and $n\in\NN$.
\end{lemma}

Using these lemmas and applying Theorems~1.9.1, 1.10.2 and 1.10.1
of~\cite{IKh}, we get Theorems~\ref{T1}, \ref{T2} and \ref{T3} respectively.

\section{Proofs of the Lemmas}

For convenience of notation, all throughout this section $C$ and $c$ denote
generic positive constants which can differ from formula to formula (and even
in the same formula), and we put $\nu=1/(p+1)$.

First of all let us fix some $\delta>0$ such that $c\,d(t)\abs t^p\leq
s(t)\leq C\,d(t)\abs t^p$ on $(-\delta {,}\delta)$, and $s(t)\geq c$ on
$[-T{,}T]\setminus(-\delta/4\,{,}\,\delta/4)$. To do so, we note that
$$
\abs{\psi(t)}\leq\abs{\psi(0)}+C\abs t^\varkappa=C\abs t^p\,\left(
\abs t^{\varkappa-p}+c\abs{\psi(0)}\abs t^{-p}\right)\leq
\min\{a{,}b\}\,\abs t^p/2
$$
for $t\in(-\delta {,}\delta)$, since $\varkappa-p>(1-p)/2>0$ and $\psi(0)=0$
in the case $p>0$. It follows clearly
$$
s(t)=d(t)\abs t^p+\psi(t)\geq\bigl(d(t)-\min\{a{,}b\}/2\bigr)\abs t^p\geq
d(t)\abs t^p/2
$$
and $s(t)\leq 2\,d(t)\abs t^p$. Finally, on the compact set
$[-T{,}T]\setminus(-\delta/4\,{,}\,\delta/4)$ we have $s(t)\geq c$ since the
function $s(\cdot)$ is continuous.

Now, let us fix some sequence $(A_n)$ such that $A_n\rightarrow+\infty$
sufficiently slowly. More precisely, we suppose that $A_n\,n^{-\nu}\rightarrow
0$ and we will give some additional conditions below. We split the interval
$[0,T]$ in three parts:
\begin{align*}
E_1&=\bigl\{t:\abs{t-\theta}<A_n\,n^{-\nu}\bigr\}=
(\theta-A_n\,n^{-\nu}\,{,}\,\theta +A_n\,n^{-\nu}),\\
E_2&=\bigl\{t:A_n\,n^{-\nu}<\abs{t-\theta}<\delta\bigr\}=
(\theta-\delta\,{,}\,\theta-A_n\,n^{-\nu})\!\cup\!
(\theta+A_n\,n^{-\nu}\,{,}\,\theta+\delta),\\
E_3&=\bigl\{t:\delta<\abs{t-\theta}\bigr\}=
(0\,{,}\,\theta-\delta)\cup(\theta+\delta\,{,}\,T).
\end{align*}

\bigskip
In order to prove Lemma~\ref{L1} we will only study the convergence of the
one-dimensional (the general case can be considered similarly) distributions
of the process
\begin{align*}
\ln Z_n(u)&=\sum_{i=1}^{n}\myint_{0}^{T}
\ln\frac{S_{\theta_u}(t)}{S_{\theta}(t)}\;dX_i(t)-n\myint_{0}^{T}
\left[\frac{S_{\theta_u}(t)}{S_{\theta}(t)}-1\right]S_{\theta}(t)\;dt\\
&=\sum_{i=1}^{n}\myint_{0}^{T}
f\;dX_i(t)-n\myint_{0}^{T}g\,S_{\theta}(t)\;dt,
\end{align*}
where we denote
$$
f=f(\theta, t, u, n)=\ln\frac{S_{\theta_u}(t)}{S_{\theta}(t)}
\quad\mbox{and}\quad
g=g(\theta, t, u, n)=\frac{S_{\theta_u}(t)}{S_{\theta}(t)}-1.
$$

The characteristic function of the random variable $\ln Z_n(u)$ can be written
as $\bigl($see, for example, Lemma~1.1 of~\cite{Kut}$\bigr)$
$$
\phi_n(\lambda)=\eb_\theta\exp\bigl\{i\,\lambda\,\ln Z_n(u)\bigr\}
=\exp\Biggl\{n\myint_{0}^{T}
\left[e^{i\,\lambda\,f}-1-i\,\lambda\,g\right]S_{\theta}(t)\;dt\Biggr\},
$$
and hence
\begin{equation}
\label{eqChFn}
\ln\phi_n(\lambda)=n\myint_{0}^{T}
\left[e^{i\,\lambda\,f}-1-i\,\lambda\,f\right]S_{\theta}(t)\;dt
+i\,\lambda\,n\myint_{0}^{T}[f-g]S_{\theta}(t)\;dt.
\end{equation}

To study this expression, let us at first establish the two following
properties.
\textsl{
\begin{itemize}
\item[{\bf a)}]For any fixed $u$, we have $\limn\limits g(\theta, t, u, n)=0$
uniformly in $\theta\in\mathbf{K}$ and $t\in E_2 \cup E_3$.
\item[{\bf b)}]We have
$$
\limn n\!\!\myint_{E_2\cup E_3}\!\!\!g^2\,S_{\theta}(t)\;dt=0.
$$
\end{itemize}
}

To prove~{\bf a)}, we put $y=t-\theta\in (E_2-\theta)\cup(E_3-\theta)$ and we
write
\begin{align*}
\abs{g(\theta,t,u,n)}&=\abs{\frac{s(t-\theta_u)}{s(t-\theta)}-1}\\
&=\abs{\frac{d(y\!-\!u\,n^{-\nu})\,\abs{y\!-\!u\,n^{-\nu}}^p
+\psi(y\!-\!u\,n^{-\nu})-d(y)\abs y^p-\psi(y)}{s(y)}}=\\
&=\frac{\Bigl|C\,\bigl(\abs{y\!-\!u\,n^{-\nu}}^p-\abs y^p\bigr)
+\psi(y\!-\!u\,n^{-\nu})-\psi(y)\Bigr|}{s(y)}\\
&\leq C\,\frac{\bigl|\abs{y\!-\!u\,n^{-\nu}}^p-\abs y^p\bigr|}{s(y)}
+\frac{\bigl|\psi(y\!-\!u\,n^{-\nu})-\psi(y)\bigr|}{s(y)}\\
&=M_1+M_2
\end{align*}
with evident notations.

For $y\in E_2-\theta$ we have
\begin{align*}
M_1&\leq C\,\frac{\bigl|\abs{y\!-\!u\,n^{-\nu}}^p-\abs y^p\bigr|}{c\abs y^p}
=C\,\biggl|\Bigl|1-\frac{u}{y\,n^\nu}\Bigr|^p-1\biggl|
\leq\frac{C\,\abs u}{\abs y\,n^\nu}\leq\frac{C\,\abs u}{A_n}\rightarrow0,\\
M_2&\leq\frac{C\,\abs{u\,n^{-\nu}}^\varkappa}{c}=C\,\abs u^\varkappa\,
n^{-\nu\varkappa}\rightarrow0\qquad\text{if $p<0$},\\
M_2&\leq\frac{C\,\abs{u\,n^{-\nu}}^\varkappa}{c{\abs y}^p}\leq
\frac{C\,\abs u^\varkappa\,n^{-\nu\varkappa}}{{(A_n\,n^{-\nu})}^p}
=C\,\abs u^\varkappa\,\frac{n^{-\nu(\varkappa-p)}}{A_n^p}\rightarrow0
\qquad\text{if $p>0$}.
\end{align*}

Finally for $y\in E_3-\theta$, using the H\"older continuity of $s(\cdot)$ we
have
$$
\abs{g(\theta,t,u,n)}\leq\frac{C\,\abs{u\,n^{-\nu}}^\varkappa}{c}=C\,
\abs u^\varkappa\,n^{-\nu\varkappa}\rightarrow0.
$$
So,~{\bf a)} is proved.

To prove~{\bf b)}, we first note that
\begin{align*}
n\myint_{E_3}g^2\,S_{\theta}(t)\;dt&=n\myint_{E_3}
\frac{\bigl(S_{\theta_u}(t)-S_\theta(t)\bigr)^2}{S_\theta(t)}\;dt
=n\!\!\myint_{E_3-\theta}\!\frac{\bigl(s(y\!-\!u\,n^{-\nu})-s(y)\bigr)^2}
{s(y)}\;dy\\
&\leq C\,n\myint_{E_3-\theta}\abs{u\,n^{-\nu}}^{2\varkappa}\;dy
=C\,\abs u^{2\varkappa}\,n^{-(2\nu\varkappa-1)}\rightarrow0
\end{align*}
since $2\nu\varkappa-1>0$. To conclude the proof it remains to show that
$$
n\!\!\!\myint_{\theta+A_n\,n^{-\nu}}^{\theta+\delta}\!\!\!\!
g^2\,S_{\theta}(t)\;dt+
n\!\!\!\myint_{\theta-\delta}^{\theta-A_n\,n^{-\nu}}\!\!\!\!
g^2\,S_{\theta}(t)\;dt\rightarrow 0.
$$

For the first term we have
\begin{align*}
n\!\!\!\myint_{\theta+A_n\,n^{-\nu}}^{\theta+\delta}\!\!\!\!
g^2\,S_{\theta}(t)\;dt
&=n\!\!\myint_{A_n\,n^{-\nu}}^{\delta}\!\!
\frac{\bigl(s(y\!-\!u\,n^{-\nu})-s(y)\bigr)^2}{s(y)}\;dy\\
&\leq n\!\!\myint_{A_n\,n^{-\nu}}^{\delta}\!\!
\frac{\bigl(s(y\!-\!u\,n^{-\nu})-s(y)\bigr)^2}{c\abs y^p}\;dy\\
&=C\,n\!\!\myint_{A_n\,n^{-\nu}}^{\delta}\!\!
\frac{\bigl(\abs{y\!-\!u\,n^{-\nu}}^p-\abs y^p\bigr)^2}{\abs y^p}\;dy\\
&\phanteq{}+C\,n\!\!\myint_{A_n\,n^{-\nu}}^{\delta}\!\!
\frac{\bigl(\psi(y\!-\!u\,n^{-\nu})-\psi(y)\bigr)^2}{\abs y^p}\;dy\\
&\phanteq{}+C\,n\!\!\myint_{A_n\,n^{-\nu}}^{\delta}\!\!
\frac{\bigl(\psi(y\!-\!u\,n^{-\nu})-\psi(y)\bigr)
\bigl(\abs{y\!-\!u\,n^{-\nu}}^p-\abs y^p\bigr)}{\abs y^p}\;dy\\
&=n\,J_1+n\,J_2+n\,J_3
\end{align*}
with evident notations. Further
$$
n\,J_1=C\,n\myint_{A_n}^{\delta\,n^\nu}\frac
{\bigl(\abs{z\!-\!u}^p-\abs z^p\bigr)^2\,n^{-2\nu p}}{\abs z^p\,n^{-\nu p}}
n^{-\nu}\,dz\leq C\myint_{A_n}^{+\infty}
\frac{\bigl(\abs{z\!-\!u}^p-\abs z^p\bigr)^2}{\abs z^p}\,dz\rightarrow0
$$
since
$$
\frac{\bigl(\abs{z\!-\!u}^p-\abs z^p\bigr)^2}{\abs z^p}
=\abs z^p\biggl(\abs{1-\frac uz}^p-1\biggr)^2\sim
\abs z^p\left(\frac Cz\right)^2=C\,\abs z^{p-2}
$$
and $p-2<-1$. Similarly
$$
n\,J_2\leq C\,n\!\!\myint_{A_n\,n^{-\nu}}^{\delta}\!\!
\frac{\abs{u\,n^{-\nu}}^{2\varkappa}}{\abs y^p}\;dy
\leq C\,\abs u^{2\varkappa}\,n^{-(2\nu\varkappa-1)}
\myint_{0}^{\delta}\abs y^{-p}\;dy\rightarrow0
$$
since $2\nu\varkappa-1>0$ and $-p>-1$. Finally
$$
\abs{n\,J_3}\leq C\sqrt{(n\,J_1)\,(n\,J_2)}\rightarrow0
$$
by Cauchy-Schwarz inequality, and so the first term converges to $0$.

The second term can be treated in the same way. So,~{\bf b)} is proved.

Now let us return to the study of the characteristic function
$\phi_n(\cdot)$. Using~\eqref{eqChFn} we can write
$$
\ln\phi_n=\varphi_1+\varphi_2+\varphi_3
$$
where we put
$$
\varphi_k=n\myint_{E_k}
\left[e^{i\,\lambda\,f}-1-i\,\lambda\,f\right]S_{\theta}(t)\;dt
+i\,\lambda\,n\myint_{E_k}[f-g]S_{\theta}(t)\;dt.
$$

For $\varphi_3$ we get
\begin{align*}
\varphi_3&=n\myint_{E_3}
\left[e^{i\,\lambda\,f}-1-i\,\lambda\,f\right]S_{\theta}(t)\;dt
+i\,\lambda\,n\myint_{E_3}[f-g]S_{\theta}(t)\;dt\\
&\simeq \frac12\,n\myint_{E_3}(i\,\lambda\,f)^2\,S_{\theta}(t)\;dt
-\frac12\,i\,\lambda\,n\myint_{E_3}g^2\,S_{\theta}(t)\;dt\\
&\simeq-\frac12\,\lambda^2\,n\myint_{E_3}g^2\,S_{\theta}(t)\;dt
-\frac12\,i\,\lambda\,n\myint_{E_3}g^2\,S_{\theta}(t)\;dt\rightarrow 0
\end{align*}
where the symbol ``$\simeq$'' means equality of limits

In the same way we get $\varphi_2\rightarrow0$, and it remains to study
$\varphi_1$. For this we put $y_u=y-u\,n^{-\nu}$, $\alpha(y)=d(y)\abs y^p$,
$$
\beta(y)=\frac{s(y)}{\alpha(y)}=1+\frac{\psi(y)}{d(y)\abs y^p}
$$
for $y\in[-T{,}T]\setminus\{0\}$, and $\beta(0)=1$. Note that the function
$\beta(\cdot)$ is clearly H\"older continuous of order
$$
\mu=\begin{cases}
\vphantom{|_{\big)}}\varkappa-p,&\text{if $p>0$}\\
\vphantom{|^{\big)}}\min\{\varkappa,-p\}, &\text{if $p<0$}
\end{cases}\quad.
$$

We have
\begin{align*}
\varphi_1&=n\myint_{E_1}
\left[e^{i\,\lambda\,f}-1-i\,\lambda\,f\right]S_{\theta}(t)\;dt
+i\,\lambda\,n\myint_{E_1}[f-g]S_{\theta}(t)\;dt\\
&=n\myint_{-A_n\,n^{-\nu}}^{A_n\,n^{-\nu}}\left[
\left(\frac{\alpha(y_u)\,\beta(y_u)}{\alpha(y)\,\beta(y)}\right)^{i\,\lambda}
-1-i\,\lambda\,\ln\frac{\alpha(y_u)\,\beta(y_u)}{\alpha(y)\,\beta(y)}
\right]\,\alpha(y)\,\beta(y)\;dy\\
&\phanteq{}-i\,\lambda\,n\myint_{-A_n\,n^{-\nu}}^{A_n\,n^{-\nu}}\left[
\frac{\alpha(y_u)\,\beta(y_u)}{\alpha(y)\,\beta(y)}
-1-\ln\frac{\alpha(y_u)\,\beta(y_u)}{\alpha(y)\,\beta(y)}
\right]\,\alpha(y)\,\beta(y)\;dy\\
&\simeq n\myint_{-A_n\,n^{-\nu}}^{A_n\,n^{-\nu}}\left[
\frac{\alpha^{i\,\lambda}(y_u)}{\alpha^{i\,\lambda}(y)}
-1-i\,\lambda\,\ln\frac{\alpha(y_u)}{\alpha(y)}
\right]\,\alpha(y)\;dy\\
&\phanteq{}+n\myint_{-A_n\,n^{-\nu}}^{A_n\,n^{-\nu}}
\frac{\alpha^{i\,\lambda}(y_u)}{\alpha^{i\,\lambda}(y)}
\left(\frac{\beta^{i\,\lambda}(y_u)}{\beta^{i\,\lambda}(y)}-1\right)
\alpha(y)\;dy\\
&\phanteq{}-i\,\lambda\,n\myint_{-A_n\,n^{-\nu}}^{A_n\,n^{-\nu}}\left[
\frac{\alpha(y_u)}{\alpha(y)}-1-\ln\frac{\alpha(y_u)}{\alpha(y)}
\right]\,\alpha(y)\;dy\\
&\phanteq{}-i\,\lambda\,n\myint_{-A_n\,n^{-\nu}}^{A_n\,n^{-\nu}}
\alpha(y_u)\left(\frac{\beta(y_u)}{\beta(y)}-1\right)\;dy\\
&=n\,I_1+n\,I_2-i\,\lambda\,n\,I_3-i\,\lambda\,n\,I_4
\end{align*}
with evident notations.

Using the H\"older continuity of $\beta(\cdot)$, we get
\begin{align*}
\abs{n\,I_4}&\leq n\myint_{-A_n\,n^{-\nu}}^{A_n\,n^{-\nu}}
\alpha(y_u)\,\frac{\abs{\beta(y_u)-\beta(y)}}{\beta(y)}\;dy
\leq n\,C\,\abs{u\,n^{-\nu}}^\mu\myint_{-A_n\,n^{-\nu}}^{A_n\,n^{-\nu}}
\frac{\alpha(y_u)}{\beta(y)}\;dy\\
&\simeq C\,\abs u^\mu\,n^{1-\nu\mu}\myint_{-A_n\,n^{-\nu}}^{A_n\,n^{-\nu}}
\alpha(y_u)\;dy
=C\,\abs u^\mu\,n^{1-\nu\mu}\myint_{(-A_n-u)\,n^{-\nu}}^{(A_n-u)\,n^{-\nu}}
d(x)\abs x^p\;dx\\
&=C\,\abs u^\mu\,n^{1-\nu\mu}\left[
\frac a{p+1}(A_n+u)^{p+1}+\frac b{p+1}(A_n-u)^{p+1}\right]
n^{-\nu(p+1)}\\
&\leq C\,\abs u^\mu\,\bigl(A_n+\abs u\bigr)^{p+1} n^{-\nu\mu}\rightarrow0
\end{align*}
if $(A_n)$ is chosen so that $A_n\,n^{-\nu^2\mu}\rightarrow0$.

Similarly, noting that $\beta^{i\,\lambda}(\cdot)$ is also H\"older continuous
of order $\mu$, we get
\begin{align*}
\abs{n\,I_2}&\leq n\myint_{-A_n\,n^{-\nu}}^{A_n\,n^{-\nu}}
\frac{\abs{\alpha^{i\,\lambda}(y_u)}}{\abs{\alpha^{i\,\lambda}(y)}}
\,\frac{\abs{\beta^{i\,\lambda}(y_u)-\beta^{i\,\lambda}(y)}}
{\abs{\beta^{i\,\lambda}(y)}}\,\alpha(y)\;dy\\
&=n\myint_{-A_n\,n^{-\nu}}^{A_n\,n^{-\nu}}
\abs{\beta^{i\,\lambda}(y_u)-\beta^{i\,\lambda}(y)}\,\alpha(y)\;dy
\leq n\,C\,\abs{u\,n^{-\nu}}^\mu\myint_{-A_n\,n^{-\nu}}^{A_n\,n^{-\nu}}
\alpha(y)\;dy\\
&=C\,\abs u^\mu\,n^{1-\nu\mu}\,A_n^{p+1}\,n^{-\nu(p+1)}=
C\,\abs u^\mu\,A_n^{p+1}\,n^{-\nu\mu}\rightarrow0
\end{align*}
under the same condition on the choice of $(A_n)$.

For $n\,I_3$ we can write
\begin{align*}
n\,I_3&=n\myint_{-A_n\,n^{-\nu}}^{A_n\,n^{-\nu}}
\left[\frac{d(y_u)\abs{y_u}^p}{d(y)\abs y^p}-1-
\ln\frac{d(y_u)\abs{y_u}^p}{d(y)\abs y^p}\right]d(y)\abs y^p\;dy\\
&=n\myint_{-A_n}^{A_n}\left[\frac{d(z\!-\!u)}{d(z)}\abs{1\!-\!\frac uz}^p-1-
\ln\Bigl(\frac{d(z\!-\!u)}{d(z)}\abs{1\!-\!\frac uz}^p\Bigr)\right]
\frac{d(z)\abs z^p}{n^{\nu(p+1)}}\;dz\\
&\rightarrow\myint_{-\infty}^{\infty}
\left[\frac{d(z\!-\!u)}{d(z)}\abs{1\!-\!\frac uz}^p-1-
\ln\Bigl(\frac{d(z\!-\!u)}{d(z)}\abs{1\!-\!\frac uz}^p\Bigr)\right]
d(z)\abs z^p\;dz.
\end{align*}
Note that the last integral is finite, since
$$
\frac{d(z\!-\!u)}{d(z)}=
\left(\frac ab\right)^{\sign(u)\,\1_{[u^-{,}u^+]}(z)}\;,
$$
and hence the integrand behaves as $C\abs z^{p-2}$ as $z\rightarrow\infty$.

Finally, for $n\,I_1$ we have
\begin{align*}
n\,I_1&=n\myint_{-A_n\,n^{-\nu}}^{A_n\,n^{-\nu}}
\left[\frac{d^{i\,\lambda}(y_u)\abs{y_u}^{i\,\lambda\,p}}
{d^{i\,\lambda}(y)\abs y^{i\,\lambda\,p}}-1-i\,\lambda\,
\ln\frac{d(y_u)\abs{y_u}^p}{d(y)\abs y^p}\right]d(y)\abs y^p\;dy\\
&=n\!\!\myint_{-A_n}^{A_n}\!\left[\Bigl(\frac{d(z\!-\!u)}{d(z)}
\abs{1\!-\!\frac uz}^p\Bigr)^{i\,\lambda}\!\!\!-1-i\,\lambda\,
\ln\Bigl(\frac{d(z\!-\!u)}{d(z)}\abs{1\!-\!\frac uz}^p\Bigr)\right]
\frac{d(z)\abs z^p}{n^{\nu(p+1)}}\;dz\\
&\rightarrow\myint_{-\infty}^{\infty}
\left[\Bigl(\frac{d(z\!-\!u)}{d(z)}
\abs{1\!-\!\frac uz}^p\Bigr)^{i\,\lambda}\!\!\!-1-i\,\lambda\,
\ln\Bigl(\frac{d(z\!-\!u)}{d(z)}\abs{1\!-\!\frac uz}^p\Bigr)\right]
d(z)\abs z^p\;dz
\end{align*}
where the last integral is finite as before.

So we get
\begin{align*}
\ln\phi_n&\rightarrow{\mathcal L}=\myint_{-\infty}^{\infty}
\left[\Bigl(\frac{d(z\!-\!u)}{d(z)}\Bigr)^{i\,\lambda}
\abs{1\!-\!\frac uz}^{i\,\lambda\,p}\!\!\!-1-i\,\lambda\,p\,
\ln\abs{1\!-\!\frac uz}\right]d(z)\abs z^p\;dz\\
&\phanteq{}\phantom{\mathcal L=}-i\,\lambda\myint_{-\infty}^{\infty}
\left[\frac{d(z\!-\!u)}{d(z)}\abs{1\!-\!\frac uz}^p-1-
p\,\ln\abs{1\!-\!\frac uz}\right]d(z)\abs z^p\;dz.
\end{align*}

To terminate the proof of Lemma~\ref{L1} it remains to show that $\mathcal
L=\ln\phi$, where $\phi(\cdot)$ is the characteristic function of $\ln Z(u)$.

Recall that
\begin{align*}
\ln Z(u)&=p\myint_{-\infty}^{+\infty}\ln\abs{1\!-\!\frac uz}\;\pi(dz)+\ln
\frac ab \myint_{0}^{u}\;Y(dz)\\
&\phanteq{}-\myint_{-\infty}^{+\infty}\biggl[\abs{1-\frac uz}^p-1-p\,
\ln\abs{1\!-\!\frac uz}\biggr] d(z)\abs z^p\;dz\\
&\phanteq{}-\frac{a\!-\!b}{p\!+\!1}\,\abs u^{p+1}\,\sign(u)\\
&=K_1+K_2-K_3-K_4
\end{align*}
with evident notations.

Hence
\begin{align*}
\ln\phi(\lambda)&=\ln\,\eb\exp\bigl\{i\,\lambda\,\ln Z(u)\bigr\}\\
&=\ln\,\eb\exp\bigl\{i\,\lambda\,K_1+i\,\lambda\,K_2\bigr\}
-i\,\lambda\,K_3-i\,\lambda\,K_4\\
&=\ln\,\eb\exp\Biggl\{i\,\lambda\myint_{-\infty}^{+\infty}
\biggl[p\,\ln\abs{1\!-\!\frac uz}+\ln\left(\frac ab\right)\,
\sign(u)\,\1_{[u^-{,}u^+]}(z)\biggr]\;Y(dz)\Biggr\}\\
&\phanteq{}+i\,\lambda\,\ln\frac ab\myint_{0}^{u}d(z)\abs z^p\;dz
-i\,\lambda\,K_3-i\,\lambda\,K_4\\
&=\myint_{-\infty}^{+\infty}
\biggl[\exp\left\{i\,\lambda\,p\,\ln\abs{1\!-\!\frac uz}
+i\,\lambda\,\ln\left(\frac ab\right)\,\sign(u)\,\1_{[u^-{,}u^+]}(z)
\right\}-1\\
&\phanteq{}\phantom{\myint_{-\infty}^{+\infty}\biggl[}
-i\,\lambda\,p\,\ln\abs{1\!-\!\frac uz}-i\,\lambda\,\ln
\left(\frac ab\right)\,\sign(u)\,\1_{[u^-{,}u^+]}(z)\biggr]d(z)\abs z^p\;dz\\
&\phanteq{}+i\,\lambda\,\ln\frac ab\myint_{0}^{u}d(z)\abs z^p\;dz
-i\,\lambda\,K_3-i\,\lambda\,K_4\\
&=\myint_{-\infty}^{+\infty}
\biggl[\abs{1\!-\!\frac uz}^{i\,\lambda\,p}
\Bigl(\frac{d(z\!-\!u)}{d(z)}\Bigr)^{i\,\lambda}-1
-i\,\lambda\,p\,\ln\abs{1\!-\!\frac uz}\biggr]d(z)\abs z^p\;dz\\
&\phanteq{}-i\,\lambda\,K_3-i\,\lambda\,K_4\\
&=\mathcal L(\lambda)+i\,\lambda\myint_{-\infty}^{+\infty}\abs{1-\frac uz}^p
\Bigl(\frac{d(z\!-\!u)}{d(z)}-1\Bigr)\,d(z)\abs z^p\;dz-i\,\lambda\,K_4\\
&=\mathcal L(\lambda)+i\,\lambda\myint_{-\infty}^{+\infty}\abs{z-u}^p\,
(a\!-\!b)\,\sign(u)\,\1_{[u^-{,}u^+]}(z)\;dz-i\,\lambda\,K_4\\
&=\mathcal L(\lambda)+i\,\lambda\,(a\!-\!b)\myint_{0}^{u}\abs{z-u}^p\;dz
-i\,\lambda\,\frac{a\!-\!b}{p\!+\!1}\,\abs u^{p+1}\,\sign(u)
=\mathcal L(\lambda).
\end{align*}

So, the convergence of the one-dimensional distributions is proved. The case
of higher-dimensional distributions can be treated similarly. The uniformity
in $\theta$ on any compact set $\mathbf{K}\subset\Theta$ is
evident. Lemma~\ref{L1} is proved.

\bigskip
Now let us prove Lemma~\ref{L2}.  For $\abs{u_1-u_2}\geq 1$ the assertion is
evident since for all $\theta$ and $n$ we have
$$
\eb_\theta
\left|Z_n^{1/2}(u_1)-Z_n^{1/2}(u_2)\right|^2\leq 4\leq
4\,\abs{u_1-u_2}^{2p+1}.
$$

Suppose now that $\abs{u_1-u_2}\leq 1$. Denoting $\Delta=u_1-u_2$ and using
Lemma~1.5 of~\cite{Kut} we can write
\begin{align*}
\eb_\theta\left|Z_n^{1/2}(u_1)-Z_n^{1/2}(u_2)\right|^2&\leq
n\myint_{0}^{T}
\left[\sqrt{S_{\theta_{u_1}}(t)}-\sqrt{S_{\theta_{u_2}}(t)}\:\right]^2\,dt\\
&=n\myint_{0}^{T}
\left[\sqrt{s(t\!-\!\theta\!-\!u_1\,n^{-\nu})}\!-\!
\sqrt{s(t\!-\!\theta\!-\!u_2\,n^{-\nu})}\right]^2\,dt\\
&=n\myint_{-\theta-u_2\,n^{-\nu}}^{T-\theta-u_1\,n^{-\nu}}
\left[\sqrt{s(y-\Delta\,n^{-\nu})}-\sqrt{s(y)}\right]^2\,dy\\
&=n\myint_{-\theta-u_2\,n^{-\nu}}^{T-\theta-u_1\,n^{-\nu}}
\frac{\bigl[s(y-\Delta\,n^{-\nu})-s(y)\bigr]^2}
{\left[\sqrt{s(y-\Delta\,n^{-\nu})}+\sqrt{s(y)}\right]^2}\,dy\\
&=n\,I_1+n\,I_2
\end{align*}
where $I_1$ and $I_2$ are the integrals of the same function over the interval
$(-\delta/2\,{,}\,\delta/2)$ and over the set
$E=(-\theta-u_2\,n^{-\nu}\,{,}\,T-\theta-u_1\,n^{-\nu})\setminus
(-\delta/2\,{,}\,\delta/2)$ respectively.


On the set $E$ we have $\abs y\geq\delta/2$, and hence
$\abs{y-\Delta\,n^{-\nu}}\geq\delta/4$ for $n$ sufficiently large. Recall that
on the set $\{y:\abs y\geq\delta/4\}$ the function $s(\cdot)$ is separated
from zero and H\"older continuous of order $\varkappa$. So, for $n$ sufficiently
large we get
$$
n\,I_2\leq n\myint_E
\frac{\abs{\Delta\,n^{-\nu}}^{2\varkappa}}{c}\;dy\leq
C\,n\abs{\Delta\,n^{-\nu}}^{p+1}=C\,\abs{u_1-u_2}^{p+1}.
$$

Further, for the first integral we have
\begin{align*}
n\,I_1&\leq C\,n\myint_{-\delta/2}^{\delta/2}\frac
{\bigl[d(y-\Delta\,n^{-\nu})\abs{y-\Delta\,n^{-\nu}}^p-d(y)\abs y^p\bigr]^2}
{\left[\sqrt{s(y-\Delta\,n^{-\nu})}+\sqrt{s(y)}\right]^2}\;dy\\
&\phanteq{}+C\,n\myint_{-\delta/2}^{\delta/2}
\frac{\bigl[\psi(y-\Delta\,n^{-\nu})-\psi(y)\bigr]^2}{s(y)}\;dy\\
&\leq C\,n\myint_{-\delta/2}^{\delta/2}\frac
{\bigl[d(y-\Delta\,n^{-\nu})\abs{y-\Delta\,n^{-\nu}}^p-d(y)\abs y^p\bigr]^2}
{\frac 12\left[\sqrt{d(y-\Delta\,n^{-\nu})}\abs{y-\Delta\,n^{-\nu}}^{p/2}
+\sqrt{d(y)}\abs y^{p/2}\right]^2}\;dy\\
&\phanteq{}+C\,n\myint_{-\delta/2}^{\delta/2}
\frac{\abs{\Delta\,n^{-\nu}}^{2\varkappa}}
{\frac 12\,d(y)\abs y^p}\;dy\\
&\leq C\,n\myint_{-\delta/2}^{\delta/2}
\left[\sqrt{d(y-\Delta\,n^{-\nu})}\abs{y-\Delta\,n^{-\nu}}^{p/2}
-\sqrt{d(y)}\abs y^{p/2}\right]^2\;dy\\
&\phanteq{}+C\,n\,\abs{\Delta\,n^{-\nu}}^{p+1}
\myint_{-\delta/2}^{\delta/2}\frac 1{d(y)}\abs y^{-p}\;dy\\
&\leq C\,n\,\abs{\Delta\,n^{-\nu}}^{p+1}\myint_{-\infty}^{\infty}
\left[\widetilde d(z-1)\abs{z-1}^{p/2}
-\widetilde d(z)\abs z^{p/2}\right]^2\;dz+C\,\abs{\Delta}^{p+1}\\
&=C\,\abs{\Delta}^{p+1}=C\,\abs{u_1-u_2}^{p+1}.
\end{align*}
Here in the last integral we have denoted $\widetilde
d(z)=\sqrt{d(z\,\Delta)}$ and noticed that the integrand
behaves as $C\abs z^{p-2}$ as $z\rightarrow\infty$.

So, in the case $\abs{u_1-u_2}\leq 1$, for all $\theta$ and $n$ sufficiently
large we get finally
$$
\eb_\theta\left|Z_n^{1/2}(u_1)-Z_n^{1/2}(u_2)\right|^2\leq
C\,n\,I_1+C\,n\,I_2\leq C\,\abs{u_1-u_2}^{p+1}.
$$
Lemma~\ref{L2} is proved.

\bigskip
It remains to verify Lemma~\ref{L3}.  Using Lemma~1.5 of~\cite{Kut}, for any
$n$, $\theta$ and $u\in U_n$ we can write
$$
\eb_\theta Z_n^{1/2}(u)\leq
\exp\biggl\{-\frac12\,n\,F\bigl(u\,n^{-\nu}\bigr)\biggr\},
$$
where for all $u\in(\alpha-\theta\,{,}\,\beta-\theta)\subset(-T{,}T)$ we
denote
$$
F(u)=\myint_{0}^{T}
\left[\sqrt{S_{\theta+u}(t)}-\sqrt{S_{\theta}(t)}\,\right]^2\;dt.
$$

First we suppose $\abs u\leq\delta/2$. Since
$\theta\in\mathbf{K}\subset(0{,}T)$, we have
\begin{align*}
F(u)&=\myint_{0}^{T}
\left[\sqrt{s(t-\theta-u)}-\sqrt{s(t-\theta)}\,\right]^2\;dt\\
&=\myint_{-\theta}^{T-\theta}\left[\sqrt{s(y-u)}-\sqrt{s(y)}\right]^2\;dy
\geq\myint_{-\varepsilon}^{\varepsilon}
\left[\sqrt{s(y-u)}-\sqrt{s(y)}\right]^2\;dy
\end{align*}
where we can take $0<\varepsilon<\delta/2$. Hence $\abs y\leq\delta/2<\delta$
and $\abs{y-u}\leq\delta$, and so we get
\begin{align*}
F(u)&\geq\myint_{-\varepsilon}^{\varepsilon}\frac{\bigl[s(y-u)-s(y)\bigr]^2}
{\left[\sqrt{s(y-u)}+\sqrt{s(y)}\right]^2}\;dy\\
&\geq\myint_{-\varepsilon}^{\varepsilon}
\frac{\Bigl[\bigl(d(y-u)\abs{y-u}^p-d(y)\abs y^p\bigr)
+\bigl(\psi(y-u)-\psi(y)\bigr)\Bigr]^2}
{\left[\sqrt{2\,d(y-u)}\abs{y-u}^{p/2}+\sqrt{2\,d(y)}\abs y^{p/2}\right]^2}
\;dy\\
&=c\myint_{-\varepsilon}^{\varepsilon}
\left[\sqrt{d(y-u)}\abs{y-u}^{p/2}-\sqrt{d(y)}\abs y^{p/2}\right]^2\;dy\\
&\phanteq{}+c\myint_{-\varepsilon}^{\varepsilon}
\frac{\bigl[\psi(y-u)-\psi(y)\bigr]^2}
{\left[\sqrt{d(y-u)}\abs{y-u}^{p/2}+\sqrt{d(y)}\abs y^{p/2}\right]^2}\;dy\\
&\phanteq{}+c\myint_{-\varepsilon}^{\varepsilon}
\frac{\bigl(\sqrt{d(y-u)}\abs{y-u}^{p/2}-\sqrt{d(y)}\abs y^{p/2}\bigr)
\bigl(\psi(y-u)-\psi(y)\bigr)}
{\sqrt{d(y-u)}\abs{y-u}^{p/2}+\sqrt{d(y)}\abs y^{p/2}}\;dy\\
&=I_1+I_2\pm\abs{I_3}
\end{align*}
with evident notations.

For the first integral we have
$$
I_1=C\abs u^p\myint_{-\varepsilon/\abs u}^{\varepsilon/\abs u}
\left[\sqrt{d\bigl(u(z-1)\bigr)}\abs{z-1}^{p/2}-\sqrt{d(u\,z)}\abs z^{p/2}
\right]^2\;dz,
$$
and so $c\abs u^p\leq I_1 \leq C\abs u^p$ since the last integral can be
bounded from above and from below by the integral of the same function over
$\RR$ and over $(-\varepsilon/T\,{,}\,\varepsilon/T)$ respectively.

For the second integral we get
$$
I_2\leq C\abs{u}^{2\varkappa}\myint_{-\varepsilon}^{\varepsilon}
\frac{1}{\left[\sqrt{d(y)}\abs y^{p/2}\right]^2}\;dy
=C\abs{u}^{2\varkappa}.
$$

Using Cauchy-Schwarz inequality, we obtain $\abs{I_3}\leq C\sqrt{I_1\,I_2}\leq
C\abs{u}^{\varkappa+\frac{p+1}2}$ for the last integral, and finally
$$
F(u)\geq c\,\abs{u}^{p+1}-C\abs{u}^{\varkappa+\frac{p+1}2}
=c\,\abs{u}^{p+1}\Bigl(1-C\,\abs{u}^{\varkappa-\frac{p+1}2}\Bigr)
\geq c\,\abs{u}^{p+1}
$$
for $u$ sufficiently small, that is for $\abs u\leq\sigma$, where $\sigma>0$
is some fixed constant.

On the other hand, we have also
$$
\inf_{\abs u\geq\sigma}F(u)=c>0,
$$
since otherwise we should have $S_{\theta+u^*}(t)=S_{\theta}(t)$ for some
fixed $u^*$ and almost all $t\in[0{,}T]$, which is impossible. Hence, for all
$\abs u\geq\sigma$ we can write
$$
F(u)\geq c\geq c\,\frac{\abs{u}^{p+1}}{T^{p+1}}=c\,\abs{u}^{p+1}.
$$

So, for all $\theta$ and $u\in(\alpha-\theta\,{,}\,\beta-\theta)$ we have
$$
F(u)\geq c\,\abs{u}^{p+1},
$$
and hence for all $n$, $\theta$ and $u\in U_n$ we can write
$$
\eb_\theta Z_n^{1/2}(u)
\leq\exp\biggl\{-\frac12\,n\,F\bigl(u\,n^{-\nu}\bigr)\biggr\}
\leq\exp\bigl\{-c\,\abs u^{p+1}\bigr\}.
$$
Lemma~\ref{L3} is proved.

\section{Concluding remarks}

\noindent\textbf{1.}\quad
For simplicity of exposition, in this paper we considered the Bayesian
estimators and the notion of asymptotic efficiency in the case of quadratic
loss function.  In fact, the results hold for a larger class of loss functions
(see~\cite{IKh} for more details).

\bigskip\noindent\textbf{2.}\quad
Again for simplicity of exposition, we considered the case where the unknown
parameter $\theta$ is a shift parameter, that is $S_\theta(t)=s(t-\theta)$.
In fact, the results hold in a more general situation, for example when the
intensity function is strictly positive (except possibly in $\theta$) and can
be written as
$$
S_\theta(t)=d(t-\theta)\abs{t-\theta}^p+\Psi(\theta,t),
$$
where $p\in(-1{,}0)\cup(0{,}1)$, the function $d(\cdot)$ is as before, and the
function $\Psi(\theta,t)$ is continuous, and uniformly in~$t$ H\"older
continuous $\bigl($of order higher than $(p+1)/2\bigr)$ with respect
to~$\theta$. In the case $p>0$ we suppose equally that
$\Psi(\theta,\theta)=0$. It is not difficult to obtain for this case the same
results as those presented above.

\bigskip\noindent\textbf{3.}\quad
Like in Chapter~6 of~\cite{IKh}, one can consider a situation when the
intensity function has several singularities of the same order.  More
precisely, we suppose that $t_1<\cdots<t_r$ with $t_r-t_1<T$, the unknown
parameter $\theta\in\Theta=(\alpha{,}\beta)\subseteq(-t_1\,{,}\,T-t_r)$, and
the intensity function is strictly positive and can be written as
$$
S_\theta(t)=\sum_{i=1}^{r}d_i(t-\theta-t_i)\,\abs{t-\theta-t_i}^p+
\Psi(\theta,t),
$$
where $p\in(-1{,}0)\cup(0{,}1)$,
$$
d_i(x)=\begin{cases}
\vphantom{|_{\big)}}a_i,&\text{if }x<0\\
\vphantom{|^{\big)}}b_i,&\text{if }x>0
\end{cases}\quad,
$$
$a_i,b_i>0$, and the function $\Psi(\theta,t)$ is continuous, and uniformly
in~$t$ H\"older continuous $\bigl($of order higher than $(p+1)/2\bigr)$ with
respect to~$\theta$. In the case $p>0$ we suppose equally that
$\Psi(\theta,\theta+t_i)=0$. It is not difficult to obtain for this problem
the results similar to those presented above. The difference is that now one
needs to introduce the process $Z$ (and hence the random variables $\zeta$ and
$\xi$) in a slightly different manner. More precisely, for each
$i=1,\ldots,n$, one should introduce a process $Z_i$ in the same manner (but
using constants $a_i$ and $b_i$ instead of $a$ and $b$) as $Z$ was
introduced. Further one should consider the process $Z$ defined by
$$
Z(u)=\prod_{i=1}^{r}Z_i(u)
$$
where the processes $Z_i$ are independent.


\end{document}